\documentclass[11pt,centertags]{article}

\usepackage{amsmath}
\usepackage{amsthm}
\usepackage{amscd}
\usepackage{amssymb}
\usepackage{verbatim}

\title{Minimal entropy rigidity for lattices in products of rank one
symmetric spaces}
\author{Christopher Connell \thanks{Supported 
in part by an NSF postdoctoral fellowship.} and Benson Farb \thanks{Supported 
in part by NSF grant DMS 9704640 and by a Sloan foundation fellowship.} 
}
\newtheorem{theorem}{Theorem}[section]
\newtheorem{theoremA}{Theorem}[section]

\refstepcounter{theoremA}
\newtheorem*{theorem*}{Theorem}
\newtheorem{proposition}[theorem]{Proposition}
\newtheorem{lemma}[theorem]{Lemma}

\newtheorem*{corollary*}{Corollary}

\theoremstyle{definition}

\newtheorem*{definition*}{Definition}

\newtheorem*{remark}{Remark}
\newtheorem*{remarks}{Remarks}

\def\eps{\epsilon}
\def\bar{\overline}
\newcommand\D{\partial}
\newcommand{\inner}[1]{\left\langle #1 \right\rangle }
\newcommand\Hess{{\rm Hess}}
\newcommand\op{\operatorname}
\newcommand\Jac{\operatorname{Jac}}
\renewcommand\deg{{\rm deg}}
\newcommand\directsum{\oplus}

\newcommand{\R}{{\mathbf \mathbb{R}}}

\renewcommand{\tilde}{\widetilde}

\DeclareMathOperator{\rank}{rk}

\DeclareMathOperator{\tr}{Tr\,}

\DeclareMathOperator{\Vol}{Vol}
\DeclareMathOperator{\vol}{Vol}
\DeclareMathOperator{\Isom}{Isom}
\DeclareMathOperator{\ent}{ent}

\DeclareMathOperator{\Id}{Id}
\DeclareMathOperator{\meas}{{\mathcal{M}}}
\DeclareMathOperator{\bary}{bar}
\DeclareMathOperator{\supp}{supp}

\begin{document}
\maketitle

\section{Introduction}

The {\em volume entropy} $h(g)$ of a closed Riemannian $n$-manifold $(M,g)$
is defined as 
$$h(g)=\lim_{R\rightarrow \infty}\frac{1}{R} \log (\Vol(B(x,R)))$$
where $B(x,R)$ is the ball of radius $R$ around a fixed point $x$ in
the universal cover $X$. (For noncompact $M$, see Section
~\ref{sec:noncompact-case}.)  The number $h(g)$ is independent of the
choice of $x$, and equals the topological entropy of the geodesic flow
on $(M,g)$ when the curvature $K(g)$ satisfies $K(g)\leq 0$ (see
\cite{Ma}).  Note that while the volume $\Vol(M,g)$ is not invariant
under scaling the metric $g$, the {\em normalized entropy}
$$\ent(g)=h(g)^n\Vol(M,g)$$ is scale
invariant.

Besson-Courtois-Gallot \cite{BCG1} showed that, if $n\geq 3$ and $M$
admits a negatively curved, locally symmetric metric $g_0$, then 
$\ent(g)$ is minimized {\em uniquely} by 
$g_0$ in the space of all Riemannian metrics on $M$.  This striking result,
called {\em minimal entropy rigidity}, 
has a great number of corollaries, including solutions to 
long-standing problems on geodesic flows, asymptotic harmonicity,
Gromov's minvol invariant, and a new proof of Mostow
Rigidity in the rank one case (see \cite{BCG2}).  

Extending minimal entropy rigidity to all
nonpositively curved, locally symmetric manifolds $M$ has
been a well-known open problem (see, e.g., \cite{BCG2}, Open Question
5).  The case of closed manifolds locally (but not necessarily globally)
isometric to products of negatively curved locally symmetric spaces of
dimension at least $3$ was announced in \cite{BCG2} and later in
\cite{BCG3}.  

In this paper we prove minimal entropy rigidity in
this case as well as in the more general setting of complete, finite
volume manifolds.  Although we haven't seen Besson-Courtois-Gallot's
proof of this result, it is likely that our proof (in the compact case)
overlaps with theirs.  In particular, we apply the powerful method 
introduced in \cite{BCG1,BCG2}, with a few new twists (see below).

\bigskip
\noindent
{\bf Statement of result. }
While the quantity $\ent(g)$ is invariant under scaling the
metric on $M$, it is not invariant under scalings of the individual
factors of $M$.  Hence the iso-entropic inequality as in \cite{BCG1,BCG2} does
not hold as stated; one needs to find the locally symmetric metric
$g_{\min}$ on $M$ which minimizes $h(g_0)^n$ for a fixed volume
$\Vol(M,g_0)$ among all locally symmetric metrics $g_0$. This
minimization problem is easily solved by Lagrange multipliers; as we show in 
\S\ref{section:bestmetric}, such a $g_{\min}$
does indeed exist and is unique up to homothety.  
Our main result is the following.

\bigskip
\begin{theoremA}[\sc\bf Minimal Entropy Rigidity]
\label{theorem:main}
{\it Let $M$ be an $n$-manifold which admits a complete, finite-volume
Riemannian metric which is locally isometric to a product of negatively
curved (rank $1$) symmetric spaces of dimension at least $3$.  Let $g$
be any other complete, finite-volume, Riemannian metric on $M$.  If $M$
is not compact then assume that $\tilde{M}$ has {\em bounded
geometry}, i.e.\ Ricci curvature bounded above and injectivity radius
bounded below.  Then

$$\ent(g)\geq \ent(g_{\min})$$
with equality iff $g$ is homothetic to the locally symmetric metric
$g_{\min}$.}
\end{theoremA}

\medskip
\noindent
{\bf Remarks. }\begin{enumerate}
\item Theorem~\ref{theorem:main} is in fact true in more generality:
  if $(N,g)$ is any finite volume Riemannian $n$-manifold of bounded
  geometry, and if $f:N\rightarrow~M$ is any proper, coarsely
  Lipschitz map, then
  $$\ent(N,g)\geq |\deg f|\ent(M,g_{\min})$$ with equality iff $f$ is a
  homothetic Riemannian covering.  It is this more general result 
that we prove.

\item P. Verovic \cite{Ve} has shown that Theorem~\ref{theorem:main}
  no longer holds if $g$ is allowed to be a Finsler metric.  This
  behavior is different than in the rank one case.
  
\item By similar arguments to those in \cite{BCG1,BCG2},
  Theorem~\ref{theorem:main} implies strong (Mostow) rigidity for the
  corresponding locally symmetric manifolds.
  
\end{enumerate}

\bigskip
\noindent
{\bf Outline of the proof of Theorem~\ref{theorem:main}. } 

In this outline we assume $M$ and $N$ are compact.  We will outline the proof 
of Remark 1 after Theorem~\ref{theorem:main}; this implies Theorem~\ref{theorem:main} by taking 
$f:N\rightarrow M$ to be the identity.  For a moment
consider the case when the metric on $N$ is nonpositively curved. 

Endow $M$ with the unique locally symmetric metric $g_{\min}$
minimizing $\ent(g)$.  Denote by $Y$ (resp. $X$) the universal cover
of $N$ (resp. $M$).  Let $\meas(\D Y), \meas(\D X)$ denote the spaces
of atomless probability measures on the visual boundaries $\D Y,\D X$
of the universal covers $Y, X$.

Morally what we do is, following the method of
\cite{BCG2}, to define a map 
$$\tilde{F}:Y\to \meas(\D
Y)\stackrel{\phi_\ast}{\to}\meas(\D
X)\stackrel{\bary}{\to}X$$
where $\phi_\ast=\D
\tilde{f}_\ast$ is the pushforward of measures and $\bary$ is the
``barycenter of a measure'' (see \S\ref{section:barycenter}).  The
inclusion $Y\to \meas(\D Y)$, denoted
$x\mapsto \nu_x$, is given by the construction of the {\em
Patterson-Sullivan measures} $\{\nu_x\}_{x\in X}$ corresponding to
$\pi_1(N)<\Isom(Y)$ (see
\S\ref{section:ps}).  An essential feature of these constructions is
that they are all canonical, so that all of the maps are {\em
  equivariant}.  Hence $\tilde{F}$ descends to a map $F:N\to M$.

One problem with this outline is that the metric on $Y$ is typically
not nonpositively curved. In fact, if $Y$ is a finite volume
nonpositively curved manifold diffeomorphic to an irreducible higher
rank locally symmetric space then $Y$ is locally symmetric (see
Chapter 9 of \cite{Eb}).  So we must find an alternative to using the
``visual boundary'' of $Y$.  This is done by constructing smooth
measures (parameterized by $s>h(g)$) on $Y$ itself, pushing them
forward via $\tilde{f}$, and convolving with Patterson-Sullivan
measure on $X$.  The maps $\tilde{F}_s$ are then defined by taking the
barycenters of these measures. In the equality case, the maps $F_s$
limit (as $s$ tends to $h(g)$ from above) to a Lipschitz map $F$ which
turns out to be the desired locally isometric covering. This idea was
first introduced in \cite{BCG1}.

As in \cite{BCG1,BCG2}, the main step in the proof is bounding the
Jacobian $|\Jac F_s|$.  A simple degree computation finishes the proof of
the compact case. The extension to the noncompact, finite volume case is
perhaps the most technical part of the proof, and requires extending
some of the ideas of \cite{BCS} to the higher rank setting to show that
$F_s$ is proper.

\bigskip
\noindent
{\bf New features. } As noted above, our proof is an application of
the method of \cite{BCG2}.  The new features which occur in
the present case are:
\begin{itemize}
\item In nonpositive curvature, Busemann functions are convex but not
strictly convex; hence a global argument is needed to prove that
the barycenter map is well-defined.
\item As shown by Albuquerque \cite{Al}, each Patterson-Sullivan
  measure $\nu_x$ on a symmetric space $X$ is supported on a subset
  $\D_F X\subset \D X$ identified with the {\em Furstenberg boundary}.
  We make essential use of this fact in several places; in particular,
  the boundary $\D_FX$ decomposes naturally as a product when
  $\Isom(X)$ does.  The visual boundary $\D X$ does not.  The
  product decomposition of $\D_FX$ allows one to reduce the key
  estimate on $|\Jac F_s|$ to the rank one case and an algebraic lemma
  of \cite{BCG2}.  Note that in the negatively curved case of
  \cite{BCG1,BCG2}, $\D_F X=\D X$.
\item In the noncompact, finite volume case, direct geometric
  estimates on $N$ and $M$ and asymptotic analysis of the maps $F_s$
  are used to show that the maps $F_s$ are proper.
\end{itemize}

\section{The best locally symmetric metric}
\label{section:bestmetric}

In this section we find the locally symmetric metric $g_{\min}$ which
minimizes $\ent(g_0)$ over all locally symmetric metrics $g_0$.  Recall
that the entropy $h(g_0)$ of a negatively curved, locally symmetric
$n$-manifold $(M,g_0)$ with maximum sectional curvature $-1$ 
equals $n+d-2$, where $d=1,2,4,$ or $8$ according to whether
$(M,g_0)$ is a real, complex, quaternionic, or Cayley hyperbolic
space.

Now suppose that $(M,g_0)$ is a complete, finite volume Riemannian
manifold which is locally isometric to a product
$$(X_1,g_1)\times\cdots\times (X_{k},g_{k})$$
of negatively curved
symmetric spaces. After possibly scaling each factor, we may assume
that $(X_i,g_i)$ has maximum sectional curvature $-1$. In that case,
$X_i$ has entropy $h_i=n_i+d_i-2$ where $n_i$ is the dimension of
$X_i$ and $d_i=1,2,4,$ or $8$ is the dimension of the division algebra
which classifies the symmetric space $X_i$. The space of locally
symmetric metrics on $M$ is then the set of metrics
$$g_\beta=\beta_1^2g_1\times\cdots\times \beta_{k}^2g_{k}$$
where each
$\beta_i$ is a positive real number.  Note that the entropy of
$\beta_i^2 g_i$ is $\beta_i^{-1}(n_i+d_i-2)$.

We now wish to minimize $\ent(g_\beta)$ over all locally symmetric
$g_\beta$ on $M$.  This is the same as minimizing $h(g_{\beta})$ while
keeping volume fixed, i.e.\ while keeping
$\prod_{i=1}^{k}\beta_i^{n_i}=1$.

It is not hard to see that the exponential growth rate of the volume
of balls in the universal cover is given by $h_0=|l|$ where
$l=\sum_{\alpha>0} m_\alpha \alpha^*$ where $\alpha^*$ is the dual of
the positive root $\alpha$ with multiplicity $m_\alpha$ (see \cite{Kn}
for a sharp asympotic for the growth of balls). In the case of
products of rank one spaces with entropies $\beta_i^{-1}h_i$, this
becomes $h(g_\beta)=\sqrt{\sum_{i=1}^{k} \beta_i^{-2}h_i^2}$, since the
roots are all orthogonal with multiplicity $h_i$.

Now, with scalings as above, we wish to minimize 
$$h_{\beta}=\sqrt{\sum_{i=1}^{k}
  \beta_i^{-2}h_i^2}=\sqrt{\sum_{i=1}^{k}\beta_i^{-2}(n_i+d_i-2)^2}$$
subject to $\prod_{i=1}^{k}\beta_i^{n_i}=1$.  An easy computation
using Lagrange multipliers now gives that the locally symmetric metric
$g_{\min}$ which minimizes $\ent(g)$ is
$$g_{\min}=\alpha_1^2g_1\times\cdots\times\alpha_{k}^2 g_{k}$$
where 
$$\alpha_i= \frac{h_i}{\sqrt{n_i}}
\prod_{i=1}^{k}\left(\frac{\sqrt{n_i}}{h_i}\right)^{\frac{n_i}{n}}$$
In this case
$$h(g_{\min})=\sqrt{n}\prod_{i=1}^{k}
\left(\frac{h_i}{\sqrt{n_i}}\right)^{n_i/n}$$

\section{Patterson-Sullivan measures on symmetric spaces}
\label{section:ps}

In this section we briefly recall Albuquerque's theory \cite{Al} of 
Patterson-Sullivan measures in higher rank symmetric spaces.  
For background on nonpositively curved manifolds, symmetric spaces,
visual boundaries, Busemann functions, etc., we refer the reader to
\cite{BGS} and \cite{Eb}.  

\subsection{Basic properties}

Let $X$ be a Riemannian symmetric space of noncompact type.  Denote by
$\D X$ the visual boundary of $X$; that is, the set of
equivalence classes of geodesic rays in $X$, endowed with the cone
topology.  Hence $X\cup \D X$ is a compactification of $X$ which
is homeomorphic to a closed ball.  Let $\Gamma$ be a lattice in
$\Isom(X)$, so that $h(g_0)<\infty$ where $(M,g_0)$ is
$\Gamma\backslash X$ with the induced metric.

Generalizing the construction of Patterson-Sullivan, Albuquerque
constructs in \cite{Al} a family of {\em Patterson-Sullivan 
measures} on $\D X$.  This is a family of 
measures $\{\nu_x\}_{x\in X}$ on $\D X$ which provide a
particularly natural embedding of $X$ into the space of
measures on $\D X$.
\begin{proposition}
\label{proposition:properties}
The family $\{\nu_x\}$ satisfies the following properties:
\begin{enumerate} 

\item Each $\nu_x$ has no atoms.

\item The family of 
measures $\{\nu_x\}$ is $\Gamma$-equivariant:
$$\gamma_*\nu_x=\nu_{\gamma x} \mbox{\ for all\ }\gamma\in\Gamma$$

\item For all $x,y\in X$, the measure $\nu_y$ is absolutely
continuous with respect to $\nu_x$.  In fact the Radon-Nikodym derivative is
given explicitly by:

\begin{align}\label{eqn:R-N}
\frac{d\nu_x}{d\nu_y}( \xi) = e^{h(g)B(x,y, \xi)}
\end{align}

\noindent
where $B(x,y,\xi)$ is the {\em Busemann function} on $X$.  For points
$x,y\in X$ and $\xi\in\D X$, the function $B:X\times X\times\D X
\rightarrow \R$ is defined by
$$B(x,y,\xi)=\lim_{t\rightarrow \infty}d_X(y,\gamma_\xi(t))-t$$
where
$\gamma_\xi$ is the unique geodesic ray with $\gamma(0)=x$ and
$\gamma(\infty)=\xi$.

\end{enumerate}
\end{proposition}

Since $\op{Isom}(X)$ acts transitively on $X$, the first property
above implies that the $\nu_x$ are all probability measures. The
second property implies no two measures are the same. Thus the
assignment $x\mapsto \nu_x$ defines an injective map
$$\nu: X\to \meas(\D X)$$
where $\meas(\D X)$ is the space
of probability measures on $X$.  Such a mapping satisfying the above
properties is called an $h(g_0)$-{\em conformal density}.
 
\subsection{The Furstenberg boundary and Albuquerque's Theorem}

The {\em Furstenberg boundary} of a symmetric space $X$ of noncompact
type is abstractly defined to be $G/P$ where $P$ is a minimal
parabolic subgroup of the connected component $G$ of the identity of
$\Isom(X)$.

Fix once and for all a basepoint $p\in X$. This choice uniquely
determines a Cartan decomposition
$\mathfrak{g}=\mathfrak{k}\oplus\mathfrak{p}$ of the Lie algebra of
$G$ where $\mathfrak{k}$ is the Lie algebra of the isotropy subgroup
$K=\op{Stab}_G(p)$ of $p$ in $G$ and $\mathfrak{p}$ is identified with
the tangent space $T_pX$.

Let $\mathfrak{a}$ be a fixed maximal abelian subspace of
$\mathfrak{p}$. The {\em rank} of $X$, denoted $\rank(X)$, is the
dimension of $\mathfrak{a}$. If $A=\exp(\mathfrak{a})$ then $A\cdot
p$ will be a maximal flat (totally geodesically embedded Euclidean
space of maximal dimension). Recall, a vector $v\in TX$ is called a
{\em regular vector} if it is tangent to a unique maximal flat.
Otherwise it is a {\em singular vector}. A geodesic is called regular
(resp. singular) if one (and hence all) of its tangent vectors are
regular (singular).  A point $\xi\in \D X$ is regular (singular) if
any (and hence all) of the geodesics in the corresponding equivalence
class are regular (singular).

A {\em Weyl chamber} $\mathfrak{a}^+$ is a choice of connected
component of the set of regular vectors in $\mathfrak{a}$. There
corresponds to $\mathfrak{a}^+$ a choice of positive roots. Similarly
if $A^+=\exp(\mathfrak{a}^+)$ then $A^+\cdot p$ is called a Weyl
chamber of the flat $A\cdot p$. The union of all the singular
geodesics in the flat $A\cdot p$ passing through $p$ is a finite set
of hyperplanes forming the boundaries of the Weyl chambers.

In our current special case of products, $\rank(X)=k$ (recall $k$
is the number of factors) since a product of geodesics, one from each
factor, produces a maximal flat.  Also, the Weyl chambers are simply
the $2^{\rank(X)}$ orthants of $\mathfrak{a}\simeq \R^{\rank(X)}$ and
$b^+$ is the unit vector corresponding to the barycenter of the
extremal unit vectors on the interior boundary of a fixed Weyl
chamber.

The Furstenberg boundary can be identified with the orbit of $G$
acting on any regular point $v(\infty)\in \D X$, the endpoint of a
geodesic tangent to a regular vector $v$. of a Weyl chamber in a fixed
flat $\mathfrak{a}$.  This follows from the fact that the action of
any such $P$ on $\D X$ fixes some regular point.

Because of this, for symmetric spaces of higher rank, behaviour on the
visual boundary can often be aptly described by its restriction to the
Furstenberg boundary.  Here we will use only some very basic
properties of this boundary.  For more details on semisimple Lie
groups and the Furstenberg boundary, see \cite{Zi}.

Let $b$ be the sum of the dual vectors of the positive roots
corresponding to $\mathfrak{a}^+$. The vector $b^+$ is called the {\em
  algebraic centroid} of $\mathfrak{a}^+$.  Set
$b^+=b/\|b\|$. 

Define the set $\D_F X\subset \D X$ to be $\D_F X=G\cdot b^+(\infty)$.
Henceforth we will refer to the Furstenberg boundary as this specific
realization. We point out that for any lattice $\Gamma$ in $\Isom(X)$,
the induced action on the boundary is transitive only on $\D_F X$.
That is, $\overline{\Gamma\cdot b^+(\infty)}=G\cdot b^+(\infty)$ even
though for any interior point $x\in X$, $\overline{\Gamma\cdot x}=\D
X$.

Theorem 7.4 and Proposition 7.5 of \cite{Al} combine to give the
following categorization which will play a crucial role in our proof
of Theorem~\ref{theorem:main}.
\begin{theorem}[Description of $\nu_x$]
\label{theorem:patterson-sullivan}
Let $(X,g_0)$ be a symmetric space of noncompact type, and let $\Gamma$ be a
lattice in $\Isom(X)$. Then
\begin{enumerate}
\item $h(g_0)=\|b\|$,
\item $b^+(\infty)$ is a regular point, and hence $\D_F X$ is a regular set,
\item For any $x\in X$, the support $\supp(\nu_x)$ of $\nu_x$ is equal
  to $\D_F X$, and
\item $\nu_x$ is the unique probability measure invariant under the
  action on $\D_F X$ of the compact isotropy group $\op{Stab}_G(x)$ at
  $x$. In particular, $\nu_p$ is the unique $K$-invariant probability
  measure on $\D_F X$.
\end{enumerate}
\end{theorem}

Note that when $X$ is has rank one, $\D_FX=\D X$. In general
$\D_F X$ has codimension ${\rank(X)}-1$ in $\D X$.

\subsection{Limits of Patterson-Sullivan measures}
We now describe the asymptotic behaviour of the $\nu_x$ as $x$ tends
to a point in $\D X$.

For any point $\xi$ of the visual boundary, let $S_\theta$ be the set
of points $\xi\in \D_F X$ such that there is a Weyl chamber $W$ whose
closure $\D \bar{W}$ in $\D X$ contains both $\theta$ and $\xi$. Let
$K_\theta$ be the subgroup of $K$ which stabilizes $S_\theta$.
$K_\theta$ acts transitively on $S_\theta$ (see the proof below).

\begin{theorem}[Support of $\nu_x$]
\label{theorem:support}
Given any sequence $\{x_i\}$ tending to $\theta\in \D X$ in the cone
topology, the measures $\nu_{x_i}$ converge in $\meas(\D_F X)$ to the
unique $K_\theta$-invariant probability measure $\nu_{\theta}$
supported on $S_\theta$.
\end{theorem}

\begin{proof} 
  Let $x_i=g_i \cdot p$, for an appropriate sequence $g_i\in G$.
  Recall that $\nu_{x_i}=(g_i)_* \nu_p$. Then combining part (4) of
  Theorem~\ref{theorem:patterson-sullivan} with Proposition 9.43 of
  \cite{GJT} have that some subsequence of the $\nu_{x_i}$ converges
  to a $K_\theta$-invariant measure $\nu_\theta$ supported on
  $S_\theta$. 
  
  Note that in \cite{GJT}, the notation $I$ refers to a subset of a
  fundamental set of roots corresponding to the face of a Weyl chamber
  containing $\theta$ in its boundary. If $g_i\cdot p=k_i a_i\cdot p$
  converges then both $k=\lim k_i$ and $a^I=lim_i a_i^I$ exist (note
  the definition of $a^I$ in \cite{GJT}). Again in the notation of
  \cite{GJT}, $K_\theta$ is the conjugate subgroup $(ka^I) K^I
  (ka^I)^{-1}$ in $K$. Moreover, $S_\theta$ is the orbit $k a^I
  K^I\cdot b^+(\infty)$.
  
  By Corollary 9.46 and Proposition 9.45 of \cite{GJT}, any other
  convergent subsequence of the $\nu_{x_i}$ produces the same measure
  in the limit, and therefore the sequence $\nu_{x_i}$ itself
  converges to $\nu_\theta$ uniquely.
\end{proof}

In the case when $\theta$ is a regular point, the above theorem
implies that $S_\theta$ is a single point and the limit measure
$\nu_\theta$ is simply the Dirac probability measure at that point
point in $\D_FX$.

\section{The barycenter of a measure}
\label{section:barycenter}

In this section we describe the natural map which is an essential 
ingredient in the method of Besson-Courtois-Gallot.  

Let $\phi$ denote the lift to universal covers of $f$ with basepoint
$p\in Y$ (resp. $f(p)\in X$), i.e. $\phi=\tilde{f}:Y\to X$.  We will
also denote the metric and Riemannian volume form on universal cover
$Y$ by $g$ and $dg$ respectively.  Then for each $s>h(g)$ and $y\in Y$
consider the probability measure $\mu_y^s$ on $Y$ in the Lebesgue
class with density given by
$$\frac{d\mu_y^s}{dg}(z)=\frac{e^{-sd(y,z)}}{\int_{Y}
  e^{-sd(y,z)}dg}.$$
The $\mu_y^s$ are well defined by the choice of
$s$.

Consider the push-forward $\phi_*\mu_y^s$, which is a measure on $X$.
Define $\sigma_y^s$ to be the convolution of $\phi_*\mu_y^s$ with the
Patterson-Sullivan measure $\nu_z$ for the symmetric metric.

More precisely, for $U\subset\partial X$ a Borel set, define
$$\sigma_y^s(U)=\int_{X}\nu_z(U)d(\phi_*\mu_y^s)(z)$$

Since $\Vert \nu_z \Vert=1$, we have
$$\Vert\sigma_y^s\Vert=\Vert\mu_y^s\Vert=1.$$

Let
$B_0(x,\theta)=B_0(\phi(p),x,\theta)$ be the Busemann function on
$X$ with respect to the basepoint $\phi(p)$ (which we
will also denote by $p$).  For $s>h(g)$ and $x\in X, y\in Y$ define a 
function
$$\mathcal{B}_{s,y}(x)=\int_{\D X}B_0(x,\theta)d\sigma_y^s(\theta)$$

By Theorem \ref{theorem:support}, the support of $\nu_z$, hence of
$\sigma_y^s$, is all of $\D_FX$, which in turn equals the
$G$-orbit $G\cdot b^+(\infty)$, where $G=\Isom(X_0)$. Hence
$$\mathcal{B}_{s,y}(x)=\int_{\D_F
X}B_0(x,\theta)d\sigma_y^s(\theta)=\int_{G\cdot
\xi}B_0(x,\theta)d\sigma_y^s(\theta)$$

Since $X$ is nonpositively curved, the Busemann function $B_0$ is
convex on $X$.  Hence $\mathcal{B}_{s,y}$ is convex on $X$, being a
convex integral of convex functions.  While $B_0$ is strictly convex
only when $X$ is negatively curved, we have the following.

\begin{proposition}[Strict convexity of $\mathcal{B}$]
  For each fixed $y$ and $s$, the function $x\mapsto
  \mathcal{B}_{y,s}(x)$ is strictly convex, and has a unique critical
  point in $X$ which is its minimum.
\end{proposition}

\begin{proof}
  For the strict convexity, it suffices to show that given a geodesic
  segment $\gamma(t)$ between two points $\gamma(0),\gamma(1)\in X$,
  there exists some $\xi\in \D_F X$ such that function
  $B_0(\gamma(t),\xi)$ is {\em strictly} convex in $t$, and hence on
  an open positive $\sigma^s_y$-measure set around $\xi$. We know it is
  convex by the comment preceding the statement of the proposition.
  
  If $B_0(\gamma(t),\xi)$ is constant on some geodesic subsegment of
  $\gamma$ for some $\xi$, then $\gamma$ must lie in some flat $\mathcal{F}$
  such that the geodesic between $\xi\in\D \mathcal{F}$ and $\gamma$
  (which meets $\gamma$ at a right angle) also lies in $\mathcal{F}$. On the
  other hand, $\xi\in \D_F X$ is in the direction of the
  algebraic centroid in a Weyl chamber, and $\gamma$ is perpendicular
  to this direction. By the properties of the roots, $\gamma$ is a
  regular geodesic (not contained in the boundary of a Weyl chamber).
  In particular, $\gamma$ is contained in exactly one flat $\mathcal{F}$.
  Furthermore, $\D_F X\cap \D \mathcal{F}$ is a finite set (an orbit of
  the Weyl group). As a result, for almost every $\xi\in \D_F X$
  $B_0(\gamma(t),\xi)$ is {\em strictly} convex in $t$.
  
  For fixed $z\in X$, by the last property listed in Proposition
  \ref{proposition:properties}, we see that 
  $$\int_{\D_F X} B_0(x,\theta)d\nu_z(\theta)$$
  tends to $\infty$ as
  $x$ tends to any boundary point $\xi\in\D X$. Then for fixed $y$ and
  $s>h(g)$, $\mathcal{B}_{y,s}(x)$ increases to $\infty$ as $x$ tends
  to any boundary point $\xi\in\D X$. Hence it has a local minimum in
  $X$, which by strict convexity must be unique.  
\end{proof}

We call the unique critical point of $\mathcal{B}_{s,y}$ the {\em
barycenter} of the measure $\sigma_y^s$, and define a map 
$\tilde{F}_s:Y\to X$ by
$$\tilde{F}_s(y)= \text{the unique critical point of
  $\mathcal{B}_{s,y}$}$$
Since for any two points $p_1,p_2\in X$
$$B_0(p_1,x,\theta)=B_0(p_2,x,\theta)+B_0(p_1,p_2,\theta),$$
$\mathcal{B}_{s,y}$ only changes by an additive constant when we
change the basepoint of $B_0$. This change does not affect the
location of critical point of $\mathcal{B}_{y,s}$. As a result,
$\tilde{F}_s$ is independent of choice of basepoints.

The equivariance of $\phi$ and of $\{\mu_y\}$ implies that
$\tilde{F}_s$ is also equivariant. Hence $\tilde{F}_s$ descends to a
map $F_s: N\to M$. As in \cite{BCG1}, we will see that $F_s$ is $C^1$,
and will estimate its Jacobian.

\section{The Jacobian estimate}
\label{sec:Jac}
The goal of this section is to prove a sharp estimate on the magnitude 
of the Jacobian of $F_s:N\rightarrow M$. 

We obtain the differential of $F_s$ by implicit differentiation:
$$0=D_{x=F_s(y)}\mathcal{B}_{s,y}(x)=\int_{\D_F X}
d{B_0}_{(F_s(y),\theta)}(\cdot) d\sigma^s_y(\theta).$$

Hence, as two forms,
\begin{gather}
0=D_y D_{x=F_s(y)}\mathcal{B}_{s,y}(x)=\int_{\D_F X}
Dd{B_0}_{(F_s(y),\theta)}(D_y F_s(\cdot),\cdot) d\sigma^s_y(\theta)\\
 -s \int_{Y}\int_{\D_F X} d{B_0}_{(F_s(y),\theta)}(\cdot)
\left< \nabla_y d_Y(y,z),\cdot\right> d\nu_{\phi(z)}(\theta) d\mu^s_y(z)
\end{gather}

The distance function $d_Y(y,z)$ is Lipschitz and $C^1$ off of the cut
locus which has Lebesgue measure 0. It follows from the Implicit
Function Theorem (see \cite{BCG2}) that $F_s$ is $C^1$ for $s>h(g)$.
By chain rule,
$$\Jac F_s = s^n\frac{\det\left(\int_{Y}\int_{\D_F X}
    d{B_0}_{(F_s(y),\theta)}(\cdot) \left< \nabla_y
      d(y,z),\cdot\right> d\nu_{\phi(z)}(\theta) d\mu^s_y(z)
  \right)}{\det\left(\int_{\D_F X}
    Dd{B_0}_{(F_s(y),\theta)}(\cdot,\cdot)
    d\sigma^s_y(\theta)\right)}.$$

Applying H{\"o}lder's inequality to the numerator gives:
$$|\Jac F_s|\leq s^n\frac{\det\left(\int_{\D_F X}
    d{B_0}_{(F_s(y),\theta)}^2 d\sigma^s_y(\theta)\right)^{1/2}
  \det\left(\int_{Y} \left<\nabla_y d_Y(y,z),\cdot\right>^2
    d\mu^s_y(z)\right)^{1/2}}{\det\left(\int_{\D_F X}
    Dd{B_0}_{(F_s(y),\theta)}(D_y F_s(\cdot),\cdot)
    d\sigma^s_y(\theta)\right)}.$$

Using that $\tr \left<\nabla_y d_Y(y,z),\cdot\right>^2 =\left|\nabla_y
d_Y(y,z)\right|^2=1$, except possibly on a measure $0$ set, we may estimate
$$\det\left(\int_{Y} \left<\nabla_y d_Y(y,z),\cdot\right>^2 
d\mu^s_y(z)\right)^{1/2}\leq \left(\frac{1}{\sqrt{n}}\right)^n.$$

Therefore
\begin{align}
\label{eq:Jac1}
|\Jac F_s|\leq
\left(\frac{s}{\sqrt{n}}\right)^n\frac{\det\left(\int_{\D_F X}
    d{B_0}_{(F_s(y),\theta)}^2
    d\sigma^s_y(\theta)\right)^{1/2}}{\det\left(\int_{\D_F X}
    Dd{B_0}_{(F_s(y),\theta)}(D_y F_s(\cdot),\cdot)
    d\sigma^s_y(\theta)\right)}.
\end{align}

\begin{theorem}[The Jacobian Estimate]
\label{JacF}
For all $s>h(g)$ and all $y\in N$ we have:
$$\vert \Jac F_s(y)\vert\le \left(\frac s{h(g_{\min})}\right)^n$$
with
equality at any $y\in N$ if and only if $D_yF_s$ is a homothety.
\end{theorem}

As discussed in the introduction, the main idea is to use the fact
that the measures $\sigma^s_y$ on $\D X$ 
are supported on
$\D_FX$, which decomposes as a product of Furtsenberg
boundaries of the rank one factors of $X$.  This can then be used to
reduce the estimate on $|\Jac F_s(y)|$ to the rank one case and an
algebraic lemma of \cite{BCG1}.

\begin{proof}
Item (1) is clear; we prove item (2).

By the hypothesis on $X$, the group $G=\Isom(X)$ can be written as 
a product $G=G_1\times G_2\cdots \times G_{\rank(X)}$, where each $G_i\neq
\op{SL}(2,\R)$ is a simple rank one Lie group.  
Theorem \ref{theorem:support} states that there exists
$\xi\in\D X$ so that for all $y\in Y$, the measure
$\sigma^s_y$ is supported on some $G$-orbit $$G\cdot \xi=\{(G_1\times
G_2\cdots \times G_{\rank(X)})\cdot \xi\}$$

Hence 
$$\D_FX=G\cdot \xi=\D_FX_1\times\cdot\times \D_FX_{\rank(X)}$$
Since each $X_i$ has rank one, $\D_FX_i=\D X_i$ so that 
$$\D_FX=\D X_1\times\cdots\times \D X_{\rank(X)}$$

Let $B_i$ denote the Busemann function for the rank one symmetric
space $X_i$ with metric $g_i$. Then for $\theta_i\in \D
X_i\subset \D X$ and $x,y\in X_i$ we have
$B_0(x,y,\theta_i)=\alpha_i B_i(x,y,\theta_i)$. Since the factors
$X_i$ are orthogonal in $X$ with respect to the metric $g_{\min}$,
the Busemann function of $(X,g_{\min})$ with basepoint $p\in X$ at
a point $\theta=(\theta_1,\dots,\theta_{\rank(X)})\in\D_F X$ is given
by
$$B_0(x,\theta)=\sum_{i=1}^{\rank(X)} \frac{\alpha_i}{\sqrt{{\rank(X)}}} B_i(x_i,\theta_i).$$

Since $\nabla_{x}^{g_{\min}}B_i=\frac{1}{\alpha_i^2}\nabla_{x}^{g_{i}}B_i$,
we may verify that $\left|\nabla_{x}^{g_{\min}}
  B_0(x,\theta)\right|^2_{g_{\min}}=1.$ Similarly,
$$\nabla_{x}^{g_{\min}} B_0=\sum_{i=1}^{\rank(X)} \frac{1}{\alpha_i
  \sqrt{{\rank(X)}}}\nabla^{g_{i}}_{x} B_i$$ Differentiating again with respect
to a fixed orthonormal basis of $T_xX$ with respect to $g_{\min}$ we obtain,
$$\Hess^{g_{\min}}_{x} B_0=\directsum_{i=1}^{\rank(X)} 
\frac{1}{\alpha_i\sqrt{{\rank(X)}}}\Hess^{g_{i}}_x B_i$$

Since $\sigma^s_y$ is supported on $G_1\times \cdots \times
G_{\rank(X)}\cdot \xi$, we can use product coordinates on $X$ to write
the right hand side of \eqref{eq:Jac1} as,

$$\left(\frac{s}{\sqrt{n}}\right)^{n} \frac{\det\left(\int_{\D_FX}
    \left(\sum_{i=1}^{\rank(X)}
      \nabla^{g_i}_{(F_s(y),\theta_i)}B_i\right)\left(\sum_{i=1}^{\rank(X)}
      \nabla^{g_i}_{(F_s(y),\theta_i)}B_i\right)^*
    d\sigma^s_y(\theta)\right)^{\frac12}}{\det\left(\int_{\D_FX}
    \directsum_{i=1}^{\rank(X)} \Hess^{g_i}_{(F_s(y),\theta_i)}B_i
    d\sigma^s_y(\theta)\right)}$$
where the superscript $*$ means
transpose and we have cancelled the factor
$\left(\frac{1}{\alpha_i\sqrt{{\rank(X)}}}\right)^n$ from the
numerator and denominator.

Recall the following estimate for positive semi-definite block
matrices, $$\det
\begin{pmatrix} A & B \cr B^* & C \end{pmatrix}\le \det(A)\det(C)$$
Then since the matrix $$\int_{\D_FX} \left(\sum_{i=1}^{\rank(X)}
  \nabla^{g_i}_{(F_s(y),\theta_i)}B_i\right)\left(\sum_{i=1}^{\rank(X)}
  \nabla^{g_i}_{(F_s(y),\theta_i)}B_i\right)^* d\sigma^s_y(\theta)$$
is positive semi-definite, by iteratively using the above estimate (on
sub-blocks) we obtain
\begin{align*}
  \det\left(\int_{\D_FX} \left(\sum_{i=1}^{\rank(X)}
      \nabla^{g_i}_{(F_s(y),\theta_i)}B_i\right)\left(\sum_{i=1}^{\rank(X)}
      \nabla^{g_i}_{(F_s(y),\theta_i)}B_i\right)^*
    d\sigma^s_y(\theta)\right)\leq\\
  \prod_{i=1}^{\rank(X)} \det \left(\int_{\D_FX}
    (\nabla^{g_i}_{(F_s(y),\theta_i)}B_i
    )(\nabla^{g_i}_{(F_s(y),\theta_i)}B_i)^*d\sigma^s_y(\theta)\right)
\end{align*}
Also, we have  
\begin{align*}
  \det\left(\int_{\D_FX} \directsum_{i=1}^{\rank(X)}
    \Hess^{g_i}_{(F_s(y),\theta_i)}B_i
    d\sigma^s_y(\theta)\right)=\prod_{i=1}^{\rank(X)} \det
  \left(\int_{\D_FX}\Hess^{g_i}_{(F_s(y),\theta_i)}B_i
    d\sigma^s_y(\theta)\right)
\end{align*}
Hence we have 
\begin{gather*}
  |\Jac F_s(y)|\leq \left(\frac{s}{\sqrt{n}}\right)^{n} \prod_{i=1}^{\rank(X)}
  \frac{\det \left(\int_{\D_FX}
      \left(\nabla^{g_i}_{(F_s(y),\theta_i)}B_i
      \right)\left(\nabla^{g_i}_{(F_s(y),\theta_i)}
B_i\right)^*d\sigma^s_y(\theta)\right)^{1/2}}{\det
    \left(\int_{\D_FX}\Hess^{g_i}_{(F_s(y),\theta_i)}B_i
      d\sigma^s_y(\theta)\right)}
\end{gather*}

Since $B_i$ is the Busemann function of rank one symmetric space
$X_i$, it follows that, setting the tensor
$H_i=(\nabla^{g_i}_{(F_s(y),\theta)} B_i)^2$, we have
$$\Hess^{g_i}_{(F_s(y),\theta_i)} B_i=\op{I}-H_i -\sum_{k=1}^{d_i-1}
J_kH_iJ_k$$
where the $J_k$ are the matrices representing the
underlying complex structure of (the division algebra corresponding
to) the symmetric space.  Lemma 5.5 of \cite{BCG2} says that for any
$n\times n$ matrix $H$ with $\tr H=1$, the following holds:
$$\frac{\displaystyle \det H^{1/2}}{\displaystyle\det \left(\op{I}-H-
  \sum_{k=1}^{d_i-1} J_k H J_k\right)}\leq \left(\frac{\sqrt{n}}{n+d-2}
\right)^{n}$$
with equality if and only if $H=\frac{1}{n}\op{I}$.  Applying
this estimate to each term of the above product now gives
$$\frac{\det \left(\int_{\D_FX}
    \left(\nabla^{g_i}_{(F_s(y),\theta_i)}B_i
    \right)\left(\nabla^{g_i}_{(F_s(y),\theta_i)}B_i
    \right)^*d\sigma^s_y(\theta) \right)^{1/2}}{\det
  \left(\int_{\D_FX}\Hess^{g_i}_{(F_s(y),\theta_i)}B_i
    d\sigma^s_y(\theta)\right)}\leq
\left(\frac{\sqrt{n_i}}{n_i+d_i-2}\right)^{n_i}$$
with equality if
and only if $H_i=\frac{1}{n_i}\op{I}$.  Hence
$$|\Jac F_s(y)|\leq \left(\frac{s}{\sqrt{n}}\right)^{n}\prod_{i=1}^{{\rank(X)}}
\left(\frac{\sqrt{n_i}}{n_i+d_i-2}\right)^{n_i}=\frac{s^n}{h(g_{\min})^n}.$$

If equality is attained then we have equality for each term. Let
$Q_i\subset T_yY$ be the subspace which is mapped by $D_yF_s$ onto
$T_{F_s(y)}X_i$.  So 
for any $v_i\in T_{F_s(y)}X_i$ and any $u_i\in Q_i$ we have  
\begin{align*}
  \left|\left<D_yF_s(u_i),v_i\right>_{g_{i}}\right| & \leq \frac{
    \sqrt{n}\sqrt{n_i}s\|v_i\|_{g_{i}}}{n_i+d_i-2}\left(
    \int_{\D_FX}dB_{y,\phi^{-1}(\theta)}(u_i)^2
    d\sigma^s_y(\theta) \right)^{1/2}
\end{align*}

Recall that $\alpha_i= \frac{h_i}{\sqrt{n_i}}
\prod_{i=1}^{\rank(X)}\left(\frac{\sqrt{n_i}}{h_i}\right)^{\frac{n_i}{n}}.$ So
multiplying each side by $\alpha_i^2$ we write the above with respect
to the metric $g_{\min}$: 
\begin{align*}
  \left|\left<D_yF_s(u_i),v_i\right>_{g_{\min}}\right| & \leq
  \frac{\alpha_i \sqrt{n}\sqrt{n_i}s\|v_i\|_{g_{\min}}}{h_i}\left(
    \int_{\D_FX}dB_{y,\phi^{-1}(\theta)}(u_i)^2
    d\sigma^s_y(\theta) \right)^{1/2}\\
  &=\sqrt{n}\|v_i\|_{g_{\min}}\frac{s}{h(g_{\min})}\left(
    \int_{\D_FX}dB_{y,\phi^{-1}(\theta)}(u_i)^2
    d\sigma^s_y(\theta) \right)^{1/2}.
\end{align*}
Hence for all $u\in T_yY$ we obtain
$$\|D_yF_s(u)\|_{g_{\min}}\leq
\sqrt{n}\frac{s}{h(g_{\min})}\left(\int_{\D_FX}
  dB_{y,\phi^{-1}(\theta)}(u)^2 d\sigma^s_y(\theta) \right)^{1/2}.$$

Now we follow Section 5 of \cite{BCG2}. It follows that for any
orthonormal basis $\{v_i\}$ of $T_{F_s(y)}X$ we have
$$\tr  (D_yF_s^*)\circ(D_yF_s)=\sum_{i=1}^n
\left<D_yF_s(v_i),D_yF_s(v_i)\right>_{g_{min}}\leq n
\left(\frac{s}{h(g_{\min})}\right)^2.$$
Lastly it follows that 
\begin{align*}
  \left(\frac{s}{h(g_{\min})}\right)^2=|\Jac F_s|^2&=\det
  (D_yF_s)^*\circ (D_yF_s)\\
  &\leq\left(\frac{\tr (D_yF_s)^*\circ (D_yF_s)}{n}\right)\leq
  \left(\frac{s}{h(g_{\min})}\right)^2.
\end{align*}

The equality implies that 
\begin{gather*}
  (D_yF_s)^*\circ (D_yF_s)=\left(\frac{s}{h(g_{\min})}\right)^2\op{I}.
\end{gather*}
Hence $D_yF_s$ is a homothety of ratio $\frac{s}{h(g_{\min})}$ as
desired.
\end{proof}

\begin{remark}
  The equality case in item (3) will not be used in what follows, but 
  will provide insight into what happens as we take limits $s\to h(g)$.
\end{remark}

\section{Finishing the proof of Theorem~\ref{theorem:main}}
\label{section:finalproof}

We now prove the statement given in Remark 1 after Theorem~\ref{theorem:main}, from which
the theorem immediately follows. We first establish that $F_s$ is
homotopic to $f$:

\begin{proposition}
\label{prop:homotopy}
For any $s>h(g)$, the map $\Psi_s\colon [0,1]\times N\to M$ defined by
$$\Psi_s(t,y)=F_{s+\frac t{1-t}}(y)$$
is a homotopy between
$\Psi_s(0,\cdot)=F_s$ and $\Psi_s(1,\cdot)=f$.
\end{proposition}

\begin{proof}
  From its definitions, $\tilde{F}_s(y)$ is continuous in $s$ and $y$.
  Observe that for fixed $y$, $\lim_{s\to
    \infty}\sigma^s_y=\nu_{\phi(y)}$. It follows that $\lim_{s\to
    \infty}\tilde{F}_s(y)=\phi(y)$. This implies the proposition.
\end{proof}

\subsection{The compact case}

Suppose $M$ and $N$ are compact.  Since for $s>h(g)$, $F_s$ is a $C^1$
map, we may simply compute the following using elementary integration
theory:

\begin{eqnarray}
  |\deg(f)| \vol(M)&=|\deg(f)| \int_{M} dg_{\min}&=\left| \int_{N}f^*dg_{\min}
\right| \cr 
&  \leq\int_{N} \left| F_s^*dg_{\min} \right|&=\int_{N} \left|
\Jac F_s dg\right| \cr 
& =\int_{N} \left| \Jac F_s\right| dg &\leq
\left(\frac{s}{h(g_{\min})}\right)^n \vol(N)
\end{eqnarray}

Letting $s\to h(g)$ gives the inequality in
Theorem~\ref{theorem:main}. In the case when equality is achieved,
after scaling the metric $g$ by the constant
$\frac{h(g)}{h(g_{min})}$, we have $h(g)=h(g_{\min})$ and
$\Vol(N)=|\deg(f)|\,\Vol(M)$. 

We remark that in the notation used in \cite{BCG1} for the sequence of
Lemmas 7.2-9 and Proposition 8.2, the measure
$\Phi^2(y,\theta)d\theta$ is simply $\sigma_y^s$ in the rank one case
(recall the Poisson kernel is $p_0(y,\theta)=e^{-h(g)
  B_0(y,\theta)}$). 

In higher rank this equivalence still holds so long as we replace
$d\theta$ by the Patterson-Sullivan measure $\nu_p$ at the basepoint
$p$.  The proofs of Lemmas 7.2-9 and Proposition 8.2 then hold in our
more general context once we replace the visual boundary by the
Furstenberg boundary and $d\theta$ by the measure $\nu_p$ as the
analysis is the same. These show that $F_s$ converges to a
nonexpanding Lipschitz map $F$ as $s\to h(g)$. Applying Theorem C.1
from Appendix C of \cite{BCG1} we obtain Theorem~\ref{theorem:main}.

\subsection{The noncompact case}
\label{sec:noncompact-case}

We now consider the case when $M$ has finite volume but is not
compact. In this setting, it is not known whether the limit in the
definition of $h(g)$ always exists. For this reason we define the
quantity $h(g)$ to be
$$h(g)=\inf \left\{ s\geq 0 \left| \ \exists C>0 \text{ such that
      }\forall x\in X, \ \int_X e^{-s d(x,z)}
    dg(z)<C\right.\right\}.$$
In fact this agrees with the previous
definition for $h(g)$ when $M$ is compact.

If $N$ has infinite volume the main theorem is automatically satisfied
so long as $h(g)>0$.  The main difficulty is that, in order for the
proof given above to work, we need to know that $F_s$ is proper (and
thus surjective since $\op{deg}(F_s)=\op{deg}(f)\neq 0$).  For this,
we will need to prove higher rank analogs of some lemmas used in
\cite{BCS} for the rank one case. For the basics of degree theory for
proper maps between noncompact spaces, see \cite{FG}.

We will show that $F_s$ is proper by essentially showing that the
barycenter of $\sigma^s_x$ lies nearby a convex set containing large
mass for this measure. This convex set is in turn far away from
$\phi(p)$ whenever $x$ is far from $p\in Y$.  We achieve this by first
estimating the concentration of the mass of $\sigma_x$ in certain
cones which will be our convex sets. One difficulty that arises in the
higher rank is that these cones must have a certain angle when
restricted to a flat. Another difficulty is that the ends of $M$ can
have large angle at infinity. In fact our methods breakdown unless we
control the asymptotic expansion of $f$ down the ends (see Remarks~
\ref{remarks:ends}).

First, we localize the barycenter of the measure $\sigma^s_x$. Let
$v_{(x,\theta)}$ be the unit vector in $S_xX$ pointing to $\theta\in\D
X$.

\begin{lemma}
\label{barycenter}
Let $K\subset X$ and $y\in Y$ be such that $(\phi_*\mu_y^s)(K)> C$ for
some constant $1>C>\frac12$. Suppose that for all $x\in X$ there
exists $v\in S_x X$ such that for all $z\in K$:
$$\int_{\D_F X}\left< v_{(x,\theta)},v \right>d\nu_z(\theta)\ge\frac1{C}-1$$ 
Then $$x\neq \tilde{F}_s(y)$$
\end{lemma}

\begin{proof}
If $\tilde{F}_s(y)=x$ then $\nabla_x \mathcal{B}_{s,y}(x)=0$.  However,
$\nabla_x \mathcal{B}_{s,y}(x)$ may be expressed as
$$\int_{X}\int_{\D_F
X}v_{(x,\theta)}d\nu_z(\theta)d\phi_*\mu_y^s(z)$$ 
where $v_{(x,\theta)}$
is the unit vector in $S_xX$ pointing to $\theta\in\D_F X$. Then we
have
\begin{align*}
  \left\Vert D_x\mathcal{B}_{s,y}\right\Vert
  &=\left\Vert\int_{X}\int_{\D_F
      X}  v_{(x,\theta)}d\nu_z(\theta)d\phi_*\mu_y^s(z) \right\Vert \\
  &\ge \left\Vert\int_K\int_{\D_F X}v_{(x,\theta)}d\nu_z
    (\theta)d\phi_*\mu_y^s(z)\right\Vert- \\
  &\hspace{2cm} \left\Vert\int_{X-K}\int_{\D_F X}
    v_{(x,\theta)}d\nu_z(\theta)d\phi_*\mu_y^s(z)\right\Vert \\
  &\ge \int_K\int_{\D_F X}\inner{v_{(x,\theta)},v}d\nu_z
  (\theta)d\phi_*\mu_y^s(z)-\phi_*\mu_y^s(X-K)\\
  &\ge \phi_*\mu_y^s(K)\left(\frac1{C}-1\right)-1+\phi_*\mu_y^s(K)\\
  &> C\left(\frac1{C}-1\right)-1+C=0\\
\end{align*}
The strictness of the inequality finishes the proof.
\end{proof}

For $v\in S X$ and $\alpha>0$ consider the convex cone,
$$E_{(v,\alpha)}=\exp_{\pi(v)}\left\{w\in T_{\pi(v)}X\ \vert\ 
  \angle_{\pi(v)} (v(\infty),w(\infty))\leq \alpha \right\},$$
where $\pi:TX \to X$ is the tangent bundle projection.

Denote by $\D E_{(v,\alpha)}\subset\D X$ its boundary at infinity.

\begin{lemma}
\label{onethird}
There exists $T_0>0$ and $\alpha_0>0$ such that for all $t\geq T_0$,
all $x\in X$, all $v\in S_xX$ and all $z\in E_{(g^{t}v,\alpha_0)}$,
$$\int_{\D_F X}\inner{v_{(x,\theta)},v}d\nu_z(\theta)\ge
\frac{\sqrt{2}}{3}.$$
\end{lemma}
\begin{proof}
  Since the isometry group of the symmetric space $X$ is transitive
  on $X$ and for any isometry $\psi$,
  $d\psi(E_{(v,\alpha)})=E_{(d\psi(v),\alpha)}$, it is sufficient to
  prove the lemma for a fixed $x$ and all $v\in S_xX$.
  
  For now choose $\alpha_0<\pi/4$. Take a monotone sequence $t_i\to
  \infty$, and any choice $z_i\in E_{(g^{t_i}v,\alpha)}$ for each
  $t_i$. It follows that some subsequence of the $z_i$, which we again
  denote by $\{z_i\}$, must tend to some point $\theta\in \D
  E_{(v,\alpha)}$.
  
  Let $\nu_\theta$ be the weak limit of the measures $\nu_{z_i}$.
  From Theorem~\ref{theorem:support}, $\nu_\theta$ is a probability
  measure supported on a set $S_\theta$ satisfying
$$\angle_x(\theta,\xi) \leq \frac{\pi}4 \quad \forall \xi\in
S_\theta.$$

Therefore we have,
\begin{align}
\label{eq:inner-product}
\int_{S_{\theta}}\inner{v_{(x,\xi)},v_{(x,\theta)}}d\nu_\theta(\xi)\ge\frac{\sqrt{2}}{2}
\end{align}  

Now whenever $\theta\in \D E_{(v,\alpha)}$ then $v=v_{(x,\theta)}+\eps
v'$ for some unit vector $v'$ and $\eps\leq \sin(\alpha)$. Using
either case above we may write
$$\int_{\D_F X}\inner{v_{(x,\xi)},v}d\nu_\theta(\xi)\ge \int_{\D_F
  X}\inner{v_{(x,\xi)},v_{(x,\theta)}}d\nu_\theta(\xi)
-\sin(\alpha).$$
So choosing $\alpha$ small enough we can guarantee
that
\begin{enumerate}
\item any two Weyl chambers intersecting $E_{(g^tv,\alpha)}$ for all
  $t>0$ in the same flat must share a common face of dimension
  $\rank(M)-1$, and
\item for any $\theta\in \D E_{(v,\alpha)}$,
  $$\int_{\D_F X}\inner{v_{(x,\xi)},v}d\nu_\theta(\xi)\ge
  \frac{\sqrt{2}}{2.5}.$$
\end{enumerate}

Let
$$E_{(v(\infty),\alpha)}=\cap_{t>0}\D E_{(g^tv,\alpha)}.$$
By the first
property used in the choice of $\alpha$ above, for any two points
$\theta_1,\theta_2\in E_{(v(\infty),\alpha)}$, either $\theta_1$ and
$\theta_2$ are in the boundary of the same Weyl chamber, or else there
is another point $\theta'$ in the intersection of the boundaries at
infinity of the closures of the respective Weyl chambers. 

By maximality there is some $\theta_0\in E_{(v(\infty),\alpha)}$
intersecting the boundary at infinity of the closure of every Weyl
chamber which intersects $E_{(g^tv,\alpha)}$ for all $t>0$.  Hence,
for every $\theta \in E_{(v(\infty),\alpha)}$, the support of the
limit measure $\nu_\theta$ satisfies $S_\theta\subset S_{\theta_0}$.
(While $\theta_0$ is not necessarily unique, the support
$S_{\theta_0}$ of the corresponding limit measure $\nu_{\theta_0}$
is.)

As $t$ increases, for any $z\in E_{(g^tv,\alpha)}$, the measures $\nu_z$
uniformly become increasingly concentrated on $S_{\theta_0}$. Then
applying the estimate \eqref{eq:inner-product} to $\theta=\theta_0$,
we may choose $T_0$ sufficiently large so that for all $z\in
E_{(g^tv,\alpha)}$ with $t>T_0$,
$$\int_{\D_F X}\inner{v_{(x,\xi)},v}d\nu_z(\xi)\ge
\frac{\sqrt{2}}{3}.$$
\end{proof}

\begin{proposition}
$F_s$ is proper.
\end{proposition}

\begin{proof}
  By way of contradiction, let $y_i\in Y$ be an unbounded sequence
  such that $\{\tilde{F}_s(y_i)\}$ lies in a compact set $K$. We may
  pass to an unbounded subsequence of $\{y_i\}$, which we again denote
  as $\{y_i\}$, such that the sequence $\phi(y_i)$ converges
  within a fundamental domain for $\pi_1(M)$ in $X$ to a point
  $\theta_0\in \D X$. Since $K$ is compact, the set
  $$A=\bigcap_{x\in K}E_{(g^{T_0}v_{(x,\theta_0)},\alpha_0)}$$
  contains an open neighborhood of $\theta_0$ and $d_X(A,K)\geq T_0$.
  Notice that $A$ is itself a cone, being the intersection of cones on
  a nonempty subset of $\D X$.
  
  We now show that $A$ contains the image $\phi(B(y_i,R_i))$ of
  increasingly large balls ($R_i\to\infty$). However, we observe from
  the fact that $A$ is a cone on an open neighborhood of $\theta_0$ in
  $\D X$ that $A$ contains balls $B(\phi(y_i),r_i)$ with
  $r_i\to\infty$. By assumption $f$, and hence $\phi$, is coarsely
  Lipschitz:
  $$d_X(\phi x,\phi y)\leq K d_Y(x,y)+C$$
  for some constants $C>0$ and
  $K\geq 1$. Therefore $\phi^{-1}(B(\phi(y_i),r_i))\supset B(y_i,R_i)$
  where $K R_i+C> r_i$. In particular $R_i\to\infty$.
  
  Hence, there exists an unbounded sequence $R_i$ such that
  $B(y_i,R_i)\subset \phi^{-1}(A)$. Furthermore, since the Ricci
  curvature is assumed to be bounded from above and the injectivity
  radius from below, we have that $\vol (B(y_i,\op{injrad}))$ is
  greater than some constant independent of $y_i$ and hence $\int_Y
  e^{-s d(y_i,z)}dg(z)>Q$ for some constant $Q>0$. By choice of $s$
  there is a constant $C_s$ depending only on $s$ such that  $\int_Y
  e^{-s d(y,z)}dg(z)<C_s$ for all $y\in Y$.

  In polar coordinates we may write,
\begin{align*}
\int_Y e^{-s d(y,z)}dg(z)&=\int_0^\infty e^{-s t}\vol(S(y,t))dt\\
&=\int_0^\infty e^{-s t}\frac{d}{dt}\vol(B(y,t))dt\\
&=-\int_0^\infty \frac{d}{dt}\left(e^{-s t}\right)\ \vol(B(y,t))dt\\
&=s\int_0^\infty e^{-s t}\vol(B(y,t))dt.
\end{align*}
Using this we may estimate, using any $\delta<s-h(g)$,
\begin{align*}
  \mu_{y_i}^s(\phi^{-1}(A))&> \mu_{y_i}^s(B(y_i,R_i))\\
  &=1-\frac{\int_{R_i}^\infty e^{-s t}\vol(B(y_i,t))dt}{\int_0^\infty
    e^{-s
      t}\vol(B(y_i,t))dt}\\
  &\geq 1-\frac{e^{-\delta R_i}\int_{R_i}^\infty e^{-(s-\delta)
      t}\vol(B(y_i,t))dt}{\int_0^\infty e^{-s
      t}\vol(B(y_i,t))dt}\\
  &\geq 1-e^{-\delta R_i}\frac{C_{s-\delta}}{Q}.
\end{align*}

Therefore for all sufficiently large $i$,
$$\mu_{y_i}^s(\phi^{-1}(A))> \frac{3}{3+\sqrt{2}}.$$
The constant
$\frac{3}{3+\sqrt{2}}$ is the constant $C$ from Lemma \ref{barycenter}
such that $\frac{1}{C}-1=\frac{\sqrt{2}}{3}$.
  
  Set $v_i=g^{T_0+1}v_{(\tilde{F}_s(y_i),\theta_0)}$. Recalling that
  $A\subset E_{(v_i, \alpha_0)}$ for all $i$, we have that for
  sufficiently large $i$,
  $$\phi_*\mu_{y_i}^s(E_{(v_i, \alpha_0)})> \frac{3}{3+\sqrt{2}}$$
  but
  $d_X(\tilde{F}_s(y_i),E_{(v_i,\alpha_0)})>T_0$, contradicting the
  conclusion of Lemma \ref{barycenter} in light of Lemma
  \ref{onethird}.
\end{proof}

\begin{remarks}
\label{remarks:ends}
\indent
\begin{enumerate}
\item In the proof of the above proposition, we used that
  $\op{injrad}$ is bounded from below and $\op{Ricci}$ curvature is
  bounded from above only to show that the volume of balls of any
  fixed radius are bounded from below. Hence this latter hypothesis
  can be replaced for the former in Theorem~\ref{theorem:main}.
  
\item Ideas from coarse topology can be used to remove the coarse
  lipschitz assumption on $f$ in the case that the ends of $M$ have
  angle at infinity bounded away from $\pi/2$. Unfortunately, $M$ may
  have ends containing pieces of flats with wide angle (consider the
  product of two rank one manifolds each with multiple cusps, or for a
  complete classification of higher rank ends see \cite{Hat}). For
  such spaces it is possible to construct a proper map $f:M\to M$ such
  that for a radial sequence $y_i\to \infty$, $\phi$ maps the bulk of
  the mass of $\mu^s_{y_i}$ into a set (almost) symmetrically arranged
  about the point $p\in X$ thus keeping $\tilde{F}_s(y_i)$ bounded.
  This explains the need for a condition on $f$ akin to the coarse
  lipschitz hypothesis.
\end{enumerate}
\end{remarks}

The inequality in Theorem~\ref{theorem:main} now follows as in the
compact case, with $\deg(f)$ and $\deg(F_s)$ suitably interpreted.

Now we complete the proof of the rigidity when we have equality in
Theorem~\ref{theorem:main}. First we assume that
$h(g)=h(g_0)$ by scaling the metric $g$ by the constant
$\frac{h(g_0)}{h(g)}.$ We note that the proofs of the lemmas in
Section 7 of \cite{BCG1} (done for the case $f=\Id$) are identical so
long as we restrict the uniformity of Lemmas 7.5 and 7.6 to be uniform
only on compact subsets. These proofs go through with only minor
modification in the case that $f$ has (local) degree $\deg(f)\neq 1$; 
this is explained in Section 8.2 of \cite{BCG1}.  In this case we
obtain the general versions of Lemma 7.6 and 7.7 of \cite{BCG1},

\begin{lemma}
  There is a subsequence $s_i$ such that the maps $F_{s_i}$ converges
  uniformly on compact sets to a continuous map $F:N\to M$ such that
  $||d_yF_{s_i}||$ is uniformly bounded on compact subsets of $N$ and
  converges to $1$ almost everywhere.
\end{lemma}

The proof of Lemma 7.8 of \cite{BCG1} then goes through without modification to obtain 

\begin{lemma}
The map $F$ is Lipschitz with Lipschitz constant less than or equal to one.
\end{lemma}

Before we proceed we must show,
\begin{lemma}
The map $F$ is proper.
\end{lemma}

\begin{proof}
By the previous lemma, $F$ is a
contracting Lipschitz map.

We note that the local notion of $deg(F)$ given by $$\deg
F(x)=\sum_{y\in F^{-1}(x)} \op{sign}(\Jac F(y))$$ is well defined for
a.e. $x\in N$. Let $P\subset M$ be the set of points which have
unbounded preimage under $F$. The set $P$ is clearly closed and of
measure 0, and $F$ acts properly on $M\setminus P$. Since $F$ is
homotopic to $f$, $\deg F(x)=\deg f(x)$ for a.e. $x$ in one connected
component $U$ of $M\setminus P$, and $\deg F(x)=0$ a.e. on the other
components. We note that Lemma C.2 and C.4 in \cite{BCG1} only require
that the injectivity radius be bounded for any finite set of
points. The proof of these two lemmas show that for every $x\in U$,
$\op{card}(F^{-1}(x)\leq \deg(f)$.

Now we will show that $F$ is proper on the closure of $U$. This implies that $U=M$. 

If the map $F$ were not proper then there would be a sequence $y_i$
tending down an end of $N$ such that $F(y_i)\in U$ limits to a point
$x_0\in P$. After pasing to a subsequence (also denoted $y_i$) we may find
compact rectifiable curves of finite length which pass through all of
the $F(y_i)$ and $x_0$. For any such $c$, by continuity, the pre-image
$F^{-1}(c)$ is therefore contained in at most $\deg(f)$ curves
$\alpha_1,\ldots,\alpha_{\deg(f)}$ one of which (say
$\alpha=\alpha_1$) can be chosen to pass through the $y_i$.

By possibly slightly perturbing the points $y_i$, Fubini's theorem guarantees that we can choose a curve $c$ such that the derivatives of $F\vert_{\alpha}$ on the pre-image curves $\alpha$ are a.e. equal to one. On the other hand, curves $\alpha$ are Lipschitz since $F$ is and therefore by the
fundamental theorem of calculus it must have the same length as $c$
which is finite. This contradicts that the $y_i$ are unbounded. 
\end{proof}

To complete the equality case (riidity) we must prove the following,
\begin{proposition} Consider two $n$-dimensional complete oriented
Riemannian manifolds of finite volume, $N$ and $M$. Suppose $F:N\to M$
is a proper Lipschitz map satisfying $d_M(F(x),F(y))\le d_N(x,y)$ for
all $x,y\in M$. Then if $\vol(N)=|\deg{f}| \vol(M)$, the map $F$ is a
Riemannian covering homotopic to $f$.  
\end{proposition}

We establish this following Appendix C of \cite{BCG1} through two
lemmas.

If we set $C(x)=\op{card}\left\{F^{-1}{x}\right\}$ The proof of Lemma C.2
in
\cite{BCG1} establishes the following.

\begin{lemma}
For almost every $y\in N$, and a.e. $x\in M$, we have $C(x)=\deg(f)$ and
$D_yF$ is an
isometry between $T_yN$ and $T_{F(y)}M$.
\end{lemma}

Since $F$ is proper, the preimage set $\left\{ F^{-1}(x) \right\}$ is
compact and hence lies in a region with injectivity radius bounded
from below. The rest of the proofs of Appendix C can then be followed
verbatim to obtain

\begin{lemma}
For every $y\in N$, $C(y)=\deg(f)$ and $F$ is a local isometry.
\end{lemma}

This lemma, in particular implies that $F$ is a Riemannian covering
map. Taking limits in Proposition \ref{prop:homotopy} establishes that
$F$ is homotopic to $f$ and hence the Proposition follows. This also
completes the proof of Theorem~\ref{theorem:main}.

\providecommand{\bysame}{\leavevmode\hbox to3em{\hrulefill}\thinspace}

\noindent
Christopher Connell:\\
Dept. of Mathematics, University of Illinois at Chicago\\
Chicago, IL 60680\\
E-mail: cconnell@math.uic.edu
\medskip

\noindent
Benson Farb:\\
Dept. of Mathematics, University of Chicago\\
5734 University Ave.\\
Chicago, Il 60637\\
E-mail: farb@math.uchicago.edu

\end{document}